\documentclass[conference,letter,10pt]{IEEEtran}
\usepackage[dvips]{graphicx}
\usepackage{float}
\usepackage{times}
\usepackage{amsmath}
\usepackage{amsfonts}
\renewcommand{\baselinestretch}{1}
\IEEEoverridecommandlockouts
\begin{document}
%
\title{A General Quaternion-valued Gradient Operator and Its Applications to Computational Fluid Dynamics and Adaptive Beamforming}


\author{\IEEEauthorblockN{Mengdi Jiang, Wei Liu}
\IEEEauthorblockA{Communications Research Group\\ Department of
Electronic and Electrical Engineering\\ University of Sheffield,
UK\\\{mjiang3, w.liu\}@sheffield.ac.uk} \and \IEEEauthorblockN{Yi Li}
\IEEEauthorblockA{School of Mathematics and Statistics\\ University of Sheffield, S3 7RH
UK\\yili@sheffield.ac.uk}\thanks{Submitted to the 2014 International Conference on Digital Signal Processing (DSP2014, held in August 2014, Hong Kong)  on 29 April 2014 and accepted on 21 May 2014.}}

\maketitle

\begin{abstract}
Quaternion-valued signal processing has received increasing attention recently. One key operation involved in derivation of all kinds of adaptive algorithms is the gradient operator. Although there have been some derivations of this operator in literature with different level of details, it is still not fully clear how this operator can be derived in the most general case and how it can be applied to various signal processing problems. In this work, we will give a general derivation of the quaternion-valued gradient operator and then apply it to two different areas. One  is to combine with the classic computational fluid dynamics (CFD) approach in wind profile prediction and the other one is to apply the result to the adaptive beamforming problem for vector sensor arrays. 
\end{abstract}


\IEEEpeerreviewmaketitle

\section{Introduction}\label{sec:introduction}

Recently, quaternion-valued signal processing has been introduced to solve  problems related to three or four-dimensional
signals, such as vector-sensor array signal processing~\cite{BihanN2004,liu14e,liu14k}, and wind profile
prediction~\cite{took09a,liu13j}. In many of the cases, the traditional complex-valued adaptive filtering operation needs to be extended to the quaternion domain to derive the corresponding adaptive algorithms. One key operation involved in derivation of quaternion-valued adaptive algorithms is the gradient operator. Although there have been some derivations of this operator in literature with different level of details, it is still not fully clear how this operator can be derived in the most general case and how it can be applied to various signal processing problems.

In this work, we will give a general derivation of the quaternion-valued gradient operator and then implement this into two different applications. One case is to combine with the classic computational fluid dynamics (CFD) approach in wind profile prediction. Wind profile prediction is a classical signal prediction problem, and we can try to solve it using traditional linear and nonlinear (neural networks) prediction techniques.
On the other hand, wind/atmospheric flow analysis is also a traditional problem in CFD,
which employs conservation laws, various physical models and numerical methods
to predict wind signals.
It might be more accurate compared to other approaches, but not without disadvantages.
For example, it is time-consuming with high computational complexity,
and it also contains uncertainties/errors in initial/boundary conditions as well as models.
Therefore, we intend to combine the two approaches in a way that retains the
efficiency of the former and the accuracy of the latter.
As a preliminary study, we will apply a quaternion-valued linear predictor
to the data generated by the CFD method to show the feasibility of
the combined approach.

Another application of quaternion is the adaptive beamforming problem in vector sensor arrays. Adaptive beamforming has been studied extensively in the past for traditional sensor array systems~\cite{vantrees02a,liu10g}. With the introduction of vector sensor arrays, such as those consisting of  crossed-dipoles and tripoles, adaptive beamforming has been extended to this area too~\cite{nehorai99a,liu14e}. A reference signal based adaptive beamformer will be set up employing the derived quaternion-valued least mean square (LMS) algorithm.

This paper is organized as follows. The general quaternion-valued gradient operator is derived  and then applied to develop the quaternion-valued LMS (QLMS) as well as the augmented QLMS (AQLMS) algorithms in Section \ref{sec:QLMS and AQLMS}. Application of the algorithms to the data generated by CFD is provided in Section \ref{sec:cfd}, and their application in adaptive beamforming is studied in Section \ref{sec:vector_sensor}. Simulation results are presented in Section \ref{sec:simulations} and conclusions are drawn in Section \ref{sec:conclusions}.

\section{Derivation of a Quaternion-valued Gradient Operator and Adaptive Filtering}
\label{sec:QLMS and AQLMS}

\subsection{Differentiation with respect to a quaternion-valued vector}\label{sec:Differentiation to a quaternion-valued vector}

We first introduce the definition of differentiation with
respect to a quaternion $q$. Assume that $f(q)$ is a function of the quaternion variable $q$, expressed as
\begin{equation}
f(q)=f_{a} + if_{b} + jf_{c} + kf_{d}\;,
\end{equation}
where $f(q)$ is in general quaternion-valued.
$f(q)$, as well as its components $f_a(q)$, $f_b(q)$, $f_c(q)$, and $f_d(q)$,
can be viewed as functions of $q_a$, $q_b$, $q_c$ and $q_d$, which can be
expressed in terms of $q$ and its involutions \cite{ell2007a}:
\begin{eqnarray}
q^{i}&=&-iqi=q_{a} + q_{b}i - q_{c}j - q_{d}k\nonumber\\
q^{j}&=&-jqj=q_{a} - q_{b}i + q_{c}j - q_{d}k\nonumber\\
q^{k}&=&-kqk=q_{a} - q_{b}i - q_{c}j + q_{d}k.
\label{eq:involutions}
\end{eqnarray}
As a
consequence, we have
\begin{eqnarray}
&~&q_{a}=\frac{1}{4}(q+q^{i}+q^{j}+q^{k})\;,
q_{b}=\frac{1}{4i}(q+q^{i}-q^{j}-q^{k})\nonumber\\
&~&q_{c}=\frac{1}{4j}(q-q^{i}+q^{j}-q^{k})\;,
q_{d}=\frac{1}{4k}(q-q^{i}-q^{j}+q^{k})\nonumber\\
\end{eqnarray}
and
\begin{eqnarray}
q+iqi+jqj+kqk&=&-2q^{*}\nonumber\\
q-iqi-jqj-kqk&=&4q_{a}.
\end{eqnarray}

Given the above relations between the involutions and the real and imaginary
parts of $q$, $f(q)$ can be regarded as a function of $q$, $q^i$, $q^j$ and
$q^k$. Therefore, in what follows, we consider generally a function of $q$ and
the involutions,
i.e., $f(q,q^{i},q^{j},q^{k})$. Using the Taylor expansion of $f$, the
differential $df$ is given by
\begin{equation}\label{eq:df1}
df=\frac{\partial f}{\partial q}dq+\frac{\partial f}{\partial q^{i}}dq^{i}+\frac{\partial f}{\partial q^{j}}dq^{j}+\frac{\partial f}{\partial q^{k}}dq^{k}
\end{equation}
Note $\partial f/\partial q$ and $dq$ are both quaternion, therefore they do
not commute. On the other hand,
$df=df_a+idf_b+jdf_c+kdf_d$.
Since $df_a$ is the differential of a real function of real numbers $q_a$,
$q_b$, $q_c$ and $q_d$, we have
\begin{eqnarray}
&~&df_a(q_a,q_b,q_c,q_d)\nonumber\\
&~&=\frac{\partial f_a}{\partial q_a}dq_a+\frac{\partial f_a}{\partial q_{b}}dq_{b}+\frac{\partial f_a}{\partial q_{c}}dq_{c}+\frac{\partial f_a}{\partial q_{d}}dq_{d}\nonumber\\
&~&=\frac{\partial f_a}{\partial q_a}[\frac{1}{4}(dq+dq^{i}+dq^{j}+dq^{k})]\nonumber\\
&~&+\frac{\partial f_a}{\partial q_{b}}[\frac{1}{4i}(dq+dq^{i}-dq^{j}-dq^{k})]\nonumber\\&~&+\frac{\partial f_a}{\partial q_{c}}[\frac{1}{4j}(dq-dq^{i}+dq^{j}-dq^{k})]\nonumber\\
&~&+\frac{\partial f_a}{\partial q_{d}}[\frac{1}{4k}(dq-dq^{i}-dq^{j}+dq^{k})]\nonumber\\
&~&=\frac{1}{4}(\frac{\partial f_a}{\partial q_a}-i\frac{\partial f_a}{\partial q_{b}}-j\frac{\partial f_a}{\partial q_{c}}-k\frac{\partial f_a}{\partial q_{d}})dq\nonumber\\
&~&+\frac{1}{4}(\frac{\partial f_a}{\partial q_a}-i\frac{\partial f_a}{\partial q_{b}}+j\frac{\partial f_a}{\partial q_{c}}+k\frac{\partial f_a}{\partial q_{d}})dq^{i}\nonumber\\
&~&+\frac{1}{4}(\frac{\partial f_a}{\partial q_a}+i\frac{\partial f_a}{\partial q_{b}}-j\frac{\partial f_a}{\partial q_{c}}+k\frac{\partial f_a}{\partial q_{d}})dq^{j}\nonumber\\
&~&+\frac{1}{4}(\frac{\partial f_a}{\partial q_a}+i\frac{\partial f_a}{\partial q_{b}}+j\frac{\partial f_a}{\partial q_{c}}-k\frac{\partial f_a}{\partial q_{d}})dq^{k}
\end{eqnarray}
Similar expressions for $idf_b$, $jdf_c$, $kdf_d$ can be derived in the same
way. The sum of the four expressions give an expressions for $df$. Comparing
the resulted expression with equation (\ref{eq:df1}), we observe that the
coefficient for $dq$ should be the same, hence:
\begin{eqnarray}
\frac{\partial f}{\partial q}&=&\frac{1}{4}(\frac{\partial f_a}{\partial q_a}
  -i\frac{\partial f_a}{\partial q_{b}}-j\frac{\partial f_a}{\partial q_{c}}-k\frac{\partial f_a}{\partial q_{d}})\nonumber\\
&~&+\frac{i}{4}(\frac{\partial f_b}{\partial q_a}-i\frac{\partial f_b}{\partial q_{b}}-j\frac{\partial f_b}{\partial q_{c}}-k\frac{\partial f_b}{\partial q_{d}})\nonumber\\
&~&+\frac{j}{4}(\frac{\partial f_c}{\partial q_a}-i\frac{\partial f_c}{\partial q_{b}}-j\frac{\partial f_c}{\partial q_{c}}-k\frac{\partial f_c}{\partial q_{d}})\nonumber\\
&~&+\frac{k}{4}(\frac{\partial f_d}{\partial q_a}-i\frac{\partial f_d}{\partial q_{b}}-j\frac{\partial f_d}{\partial q_{c}}-k\frac{\partial f_d}{\partial q_{d}})\nonumber\\
&=&\frac{1}{4}(\frac{\partial f}{\partial q_a}-\frac{\partial f}{\partial q_{b}}i-\frac{\partial f}{\partial q_{c}}j-\frac{\partial f}{\partial q_{d}}k)
\end{eqnarray}
Therefore, $\dfrac{\partial f(q)}{\partial q}$ is given by
\begin{equation}
\dfrac{\partial f(q)}{\partial q}=\frac{1}{4}(\displaystyle\frac{\partial{f(q)}}{\partial q_a}-\displaystyle\frac{\partial{f(q)}}{\partial q_b} i-\displaystyle\frac{\partial{f(q)}}{\partial q_c} j- \displaystyle\frac{\partial{f(q)}}{\partial q_d} k)
\label{eq:general_definition}
\end{equation}
Expressions for $\partial f/\partial q^i$, $\partial f/\partial q^j$ and
$\partial f/\partial q^k$ can be derived similarly. Note that, in general
$\partial f(q)/ \partial q_b$ is a quaternion, therefore $\partial
f(q)/\partial q_b i \neq i \partial f(q) /\partial q_b$, i.e., the two factors
do not commute. The same argument applies to the last two terms in equation
(\ref{eq:general_definition}).

$f(q)$ can also be viewed as a function of $q^*$ and its
involutions. Following the same arguments, we can also find the
derivative of $f(q)$ with respect to $q^{*}$, which is given by
\begin{equation}
\dfrac{\partial f(q)}{\partial q^{*}}=\frac{1}{4}(\displaystyle\frac{\partial{f(q)}}{\partial q_a}+\displaystyle\frac{\partial{f(q)}}{\partial q_b} i+\displaystyle\frac{\partial{f(q)}}{\partial q_c} j+ \displaystyle\frac{\partial{f(q)}}{\partial q_d} k)
\label{eq:conj_general_definition}
\end{equation}
where $q^{*}=q_a-q_{b}i-q_{c}j-q_{d}k$.

With these results, we can then calculate the derivatives of some simple
quaternion functions. For example, we easily obtain
\begin{equation}
 \frac{\partial q}{\partial q}=1,~\frac{\partial q}{\partial q^{*}}=-\frac{1}{2}\;.
 \end{equation}
On the other hand, the product rule is not true in general due to the
non-commutativity of quaternion products. However, it holds for
the differentiation of quaternion-valued functions to real variables.
Suppose $f(q)$ and $g(q)$ are two quaternion-valued functions of the quaternion variable $q$,
and $q_a$ is the real variable. Then we can have the following result
\begin{eqnarray}
\frac{\partial f(q)g(q)}{\partial q_a}&=&\frac{\partial }{\partial q_a}(f_a+if_b+jf_c+kf_d)g\nonumber\\
&=&\frac{\partial f_a g}{\partial q_a}+i\frac{\partial f_b g}{\partial q_a}+j\frac{\partial f_c g}{\partial q_a}+k\frac{\partial f_d g}{\partial q_a}\nonumber\\
&=&(f_a\frac{\partial g}{\partial q_a}+\frac{\partial f_a}{\partial q_a}g)+i(f_b\frac{\partial g}{\partial q_a}+\frac{\partial f_b}{\partial q_a}g)\nonumber\\
&~&+j(f_c\frac{\partial g}{\partial q_a}+\frac{\partial f_c}{\partial q_a}g)+k(f_d\frac{\partial g}{\partial q_a}+\frac{\partial f_d}{\partial q_a}g)\nonumber\\
&=&(f_a+if_b+jf_c+kf_d)\frac{\partial g}{\partial q_a}\nonumber\\
&~&+(\frac{\partial f_a}{\partial q_a}+i\frac{\partial f_b}{\partial q_a}+j\frac{\partial f_c}{\partial q_a}+k\frac{\partial f_d}{\partial q_a})g\nonumber\\
&=&f(q)\frac{\partial g(q)}{\partial q_a}+\frac{\partial f(q)}{\partial q_a} g(q)
\end{eqnarray}

When the quaternion variable $q$ is replaced by a quaternion-valued vector $\textbf{w}$, given by
\begin{equation}
\textbf{w} = [w_1~w_2~\cdots~w_{M}]^{T}
\end{equation}
where $w_m = a_m+b_mi+c_mj+d_mk$, $m=1, ..., M$, the differentiation of $f(\textbf{w})$ with respect to $\textbf{w}$ can be derived using a combination of (\ref{eq:general_definition}) as follows
\begin{eqnarray}
\dfrac{\partial f}{\partial \textbf{w}}=\frac{1}{4}\left[\begin{matrix}
              \frac{\partial f}{\partial a_1}-\frac{\partial f}{\partial b_1} i-\frac{\partial f}{\partial c_1} j-\frac{\partial f}{\partial d_1} k\\
              \frac{\partial f}{\partial a_2}-i\frac{\partial f}{\partial b_2} i-\frac{\partial f}{\partial c_2} j-\frac{\partial f}{\partial d_2} k\\
              \vdots \\
              \frac{\partial f}{\partial a_{M}}-\frac{\partial f}{\partial b_{M}} i-\frac{\partial f}{\partial c_{M}} j-\frac{\partial f}{\partial d_{M}} k
             \end{matrix}\right]
\label{eq:vector_definition}
\end{eqnarray}
Similarly, we define $\dfrac{\partial f}{\partial \textbf{w}^{*}}$ as
\begin{eqnarray}
\dfrac{\partial f}{\partial \textbf{w}^{*}}=\frac{1}{4}\left[\begin{matrix}
              \frac{\partial f}{\partial a_1}+\frac{\partial f}{\partial b_1} i+\frac{\partial f}{\partial c_1} j+\frac{\partial f}{\partial d_1} k\\
              \frac{\partial f}{\partial a_2}+\frac{\partial f}{\partial b_2} i+\frac{\partial f}{\partial c_2} j+\frac{\partial f}{\partial d_2} k\\
              \vdots \\
              \frac{\partial f}{\partial a_{M}}+\frac{\partial f}{\partial b_{M}} i+\frac{\partial f}{\partial c_{M}} j+ \frac{\partial f}{\partial d_{M}} k
             \end{matrix}\right]
\label{eq:conj_vector_definition}
\end{eqnarray}
Obviously, when $M=1$, (\ref{eq:vector_definition}) and (\ref{eq:conj_vector_definition}) are reduced to (\ref{eq:general_definition}) and (\ref{eq:conj_general_definition}), respectively.

\subsection{The QLMS algorithm}
The output $y[n]$ and error $e[n]$ of a standard adaptive filter can be expressed as
\begin{eqnarray}
y[n]&=&{\textbf{w}^{T}[n]}{\textbf{x}[n]}\\
e[n]&=&d[n]-{\textbf{w}^{T}[n]}{\textbf{x}[n]},
\end{eqnarray}
where $\textbf{w}[n]$ is the adaptive weight vector with a length of $M$, $d[n]$ is the reference signal, $\textbf{x}[n]=[x[n-1], x[n-2], \cdots, x[n-M]]^{T}$ is the input sample sequence, and $\{\cdot\}^{T}$ denotes the transpose operation. 
The cost function with the quaternion-valued error
is $J_0[n]=e[n]e^{*}[n]$. Its gradient is given by
\begin{eqnarray}
\nabla_{\textbf{w}^{*}}J_0[n]=\frac{\partial {J_0[n]}}{\partial \textbf{w}^{*}}\\
\label{eq:conj_gradient_cost_function}
\nabla_{\textbf{w}}J_0[n]=\frac{\partial {J_0[n]}}{\partial \textbf{w}}
\label{eq:gradient_cost_function}
\end{eqnarray}
with respect to $\textbf{w}^{*}[n]$ and $\textbf{w}[n]$, respectively. According to \cite{mandic2011a,brandwood83a},
the conjugate gradient gives the maximum steepness direction for the optimization surface.
Therefore, the conjugate gradient $\nabla_{\textbf{w}^{*}}J_0[n]$ will be used to derive the update of the
coefficient weight vector.

First we have
\begin{eqnarray}
J_0[n]=d[n]d^{*}[n]-d[n]{\textbf{x}^{H}[n]}{\textbf{w}^{*}[n]}-{\textbf{w}^{T}[n]} {\textbf{x}[n]}d^{*}[n]\nonumber\\
+{\textbf{w}^{T}[n]} {\textbf{x}[n]}{\textbf{x}^{H}[n]}{\textbf{w}^{*}[n]}
\label{eq:extended_cost_function}
\end{eqnarray}
For different parts, we obtain the following results
\begin{equation}
\frac {\partial (d[n]d^{*}[n])}{\partial {\textbf{w}^{*}[n]}} = 0
\label{eq:part_1}
\end{equation}
\begin{equation}
\frac {\partial (d[n]{\textbf{x}^{H}[n]}{\textbf{w}^{*}[n]})}{\partial {\textbf{w}^{*}[n]}} = d[n]\textbf{x}^{*}[n]
\label{eq:part_2}
\end{equation}
\begin{equation}
\frac {\partial ({\textbf{w}^{T}[n]}{\textbf{x}[n]}d^{*}[n])}{\partial {\textbf{w}^{*}[n]}} = -\frac{1}{2}d[n]\textbf{x}^{*}[n]
\label{eq:part_3}
\end{equation}
\begin{equation}
\frac{\partial({\textbf{w}^{T}[n]} {\textbf{x}[n]}{\textbf{x}^{H}[n]}{\textbf{w}^{*}[n]})}{\partial {\textbf{w}^{*}[n]}}=\frac{1}{2}{\textbf{w}^{T}[n]}{\textbf{x}[n]}{\textbf{x}^{*}[n]}
\label{eq:part_4}
\end{equation}
Then we have the final gradient result
\begin{equation}
\nabla_{\textbf{w}^{*}}J_0[n]=-\frac{1}{2}e[n]\textbf{x}^{*}[n].
\end{equation}
With the general update equation for the weight vector
\begin{equation}
\textbf{w}[n+1] = \textbf{w}[n]-\mu \nabla_{\textbf{w}^{*}}J_0[n],
\end{equation}
we arrive at the following update equation for the QLMS algorithm with step size $\mu$
\begin{equation}
\textbf{w}[n+1] = \textbf{w}[n]+\mu(e[n]\textbf{x}^{*}[n]).
\label{eq:update_weight_vector}
\end{equation}

\subsection{The AQLMS algorithm}

Recently, to fully exploit the second-order statistics of the signals, an augmented formulation of the data vector has been proposed, first for complex-valued signals and then for quaternion-valued ones. For complex-valued signals, the augmented vector is composed of the original data and its conjugate, while for the latter, due to existence of the three perpendicular quaternion involutions, the choice for the augmented vector is not unique. Without loss of generality, here we adopt the simplest formulation  by combining the data vector $\textbf{x}[n]$ and its conjugate $\textbf{x}^{*}[n]$ to produce an augmented vector $\textbf{x}_{a}[n]=\big[\textbf{x}^{T}[n]~~\textbf{x}^{H}[n]\big]^{T}$~\cite{took10a}, where $\{\cdot\}^{H}$ is a combination of the operations of $\{\cdot\}^{T}$ and $\{\cdot\}^{*}$ for a quaternion. For such a ``widely linear'' model, the quaternion-valued output for the conjugate part of the input is given by
\begin{equation}
\hat{y}[n]=\textbf{g}^{T}[n]\textbf{x}^{*}[n],
\label{eq:aug_output}
\end{equation}
where $\textbf{g}[n]$ denotes the weight vector for the conjugate part of the input $\textbf{x}[n]$.

As to the AQLMS algorithm, the update of the weight vector of the conjugate part $\textbf{g}[n]$
can be found with the same method as that of the QLMS in (\ref{eq:update_weight_vector}), i.e.
\begin{equation}
\textbf{g}[n+1]=\textbf{g}[n]+\mu(e[n]\textbf{x}[n]).
\label{eq:conj_update_weight_vector}
\end{equation}
With the augmented weight vector $\textbf{h}_{a}[n]$ defined as
\begin{equation}
\textbf{h}_{a}[n]=\big[\textbf{w}^{T}[n]~~\textbf{g}^{T}[n]\big]^{T},
\label{eq:aug_weight}
\end{equation}
we obtain the following update equation
\begin{equation}
\textbf{h}_{a}[n+1]=\textbf{h}_{a}[n]+\mu(e_{a}[n]{\textbf{x}_{a}}^{*}[n])
\label{eq:augmentedweight}
\end{equation}
where $e_{a}[n]=d[n]-{\textbf{h}_{a}}^{T}[n]\textbf{x}_{a}[n]$.

\section{Application to CFD data}\label{sec:cfd}

CFD is a branch of fluid mechanics. It uses numerical approaches to solve fluid flow problems.

\subsection{Fluid Dynamics Equations}\label{sec:cfd_equations}

The Navier-Stokes equations are the basis of fluid problems.  
They are essentially the mathematical formulation of the Newton's second law
applied to fluid motions.
The general expression of the equations is
\begin{equation}
\rho(\frac{\partial \textbf{u}}{\partial t}+(\textbf{u} \cdot \nabla) \textbf{u})=-\nabla{P}+\eta \Delta{\textbf{u}}
\label{eq:cfd}
\end{equation}
where $\textbf{u}$ is the fluid velocity at a particular spatial location at a
given time, $P$ is the pressure and $\rho$ is the fluid density.
The left hand side of the equation is the acceleration of the fluid,
whilst on the right side are (the gradient of) the forces, including pressure and viscous force.
Together with the conservation of mass and suitable boundary conditions,
the Navier-Stokes equations can model a large class of fluid motions accurately~\cite{Ferziger2001a}.

\subsection{Turbulence}

The second term on the left hand side of equation (\ref{eq:cfd}) represents the
contribution from the advection of fluid particles to fluid acceleration,
and is customarily called
the inertial force. The second term on the right hand side represents the
viscous force.
The ratio of these two forces is defined as the Reynolds
number ($Re$).
As it turns out, when $Re$ is large, the flows tend to become unstable and
generate a spectrum of high frequency components in the velocity signal.
Such a regime of fluid motions is called turbulence.
Atmospheric flows, including
the wind fields around wind farms, are always turbulence~\cite{Pope2000a}. Due
to the presence of the high frequency components, the CFD calculation of the
velocity signal in turbulent wind fields becomes very time consuming unless
simplifying models are introduced.

\subsection{Direct numerical simulation (DNS)}

DNS solves the Navier-Stokes equations directly without any turbulence models.
The advantage of this method is that it is simple as well as accurate with complete
information. However, the computational cost can be very high if $Re$ is large.
Therefore, this method
is not yet applicable to practical situations, for example, the atmospheric
flows we will deal with \cite{Ferziger2001a}. Nevertheless, as a first step,
we choose to use DNS to generate the velocity signals in this study.

\subsection{Data Generation Using CFD}

The velocity signals are generated by DNS, where the Navier-Stokes equations are solved
using a pseudo-spectral method.
The CFD code is written in FORTRAN 90.
Running the code, we generate a time series of
three dimensional turbulent wind velocity fields in a 3-D periodic box.
We consider the flow field as an idealized wind field with the mean velocity
having been subtracted, and the signal normalized.


\section{Application to Adaptive Beamforming}\label{sec:vector_sensor}
\subsection{Quaternionic array signal model}

\begin{figure}[http]
\begin{center}
   \includegraphics[width=0.65\linewidth]{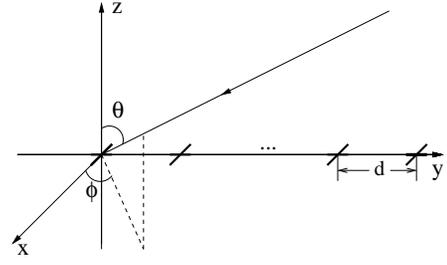}
   \caption{A ULA with crossed-dipoles.
    \label{fig:ula}}
\end{center}
\end{figure}
A uniform linear array (ULA) with $M$ crossed-dipole pairs is shown in Fig.~\ref{fig:ula}. These pairs are located along the y-axis with an adjacent spacing $d$, and at each location the two crossed components are parallel to x-axis and y-axis, respectively. Assume there is a far-field incident signal with direction of arrival (DOA) defined by the angles $\theta$ and $\phi$ impinging upon the array from the y-z plane, so that $\phi=\pi/2$ as well as $\phi=-\pi/2$, and $0 \leq \theta \leq \pi/2$. As a result, the spatial steering vector for the signal is expressed as
\begin{eqnarray}
\textbf{S}_c(\theta,\phi)&=&[1,e^{-j2\pi d sin{\theta} sin{\phi}/{\lambda}},\nonumber\\
&~&\cdots, e^{-j2\pi (M-1)d sin{\theta} sin{\phi}/{\lambda}}]^{T}
\end{eqnarray}
where $\lambda$ is the wavelength of the incident signal. For a crossed dipole the spatial-polarization coherent vector can be given by~\cite{compton81a,li91a}
\begin{equation}
\textbf{S}_p(\theta,\phi,\gamma,\eta) =
\begin{cases}
     [ -cos{\gamma},cos{\theta} sin{\gamma} e^{j\eta} ] & \text{for $\phi=\pi/2$} \\
     [ cos{\gamma},-cos{\theta} sin{\gamma} e^{j\eta} ] & \text{for $\phi=-\pi/2$}
\end{cases}
\end{equation}
where $\gamma$ is the auxiliary polarization angle with $\gamma \in [0,\pi/2]$, and the $\eta \in [-\pi,\pi]$ is the polarization phase difference.

The array structure can be divided into two sub-arrays. One is parallel to the x-axis and the other is parallel to the y-axis. The complex-valued steering vector of the x-axis sub-array is given by
\begin{equation}
\textbf{S}_x(\theta,\phi,\gamma,\eta) =
\begin{cases}
      -cos{\gamma}\textbf{S}_c(\theta,\phi) & \text{for $\phi=\pi/2$} \\
      cos{\gamma}\textbf{S}_c(\theta,\phi) & \text{for $\phi=-\pi/2$}
\end{cases}
\end{equation}
and for the y-axis it is expressed as
\begin{equation}
\textbf{S}_y(\theta,\phi,\gamma,\eta) =
\begin{cases}
      cos{\theta} sin{\gamma} e^{j\eta}\textbf{S}_c(\theta,\phi) & \text{for $\phi=\pi/2$} \\
      -cos{\theta} sin{\gamma} e^{j\eta}\textbf{S}_c(\theta,\phi) & \text{for $\phi=-\pi/2$}
\end{cases}
\end{equation}
Combining these two steering vectors together, we have a quaternion-valued composite steering vector given as below
\begin{equation}
\textbf{S}_q(\theta,\phi,\gamma,\eta)=\textbf{S}_x(\theta,\phi,\gamma,\eta)+i\textbf{S}_y(\theta,\phi,\gamma,\eta).
\end{equation}

The response of the array is
\begin{eqnarray}
r(\theta,\phi,\gamma,\eta)=\textbf{w}^{H}\textbf{S}_q(\theta,\phi,\gamma,\eta)
\end{eqnarray}
where $\textbf{w}$ is the quaternion-valued weight vector.

\subsection{Reference signal based adaptive beamforming}
\begin{figure}[htbp]
\begin{center}
   \includegraphics[width=0.65\linewidth]{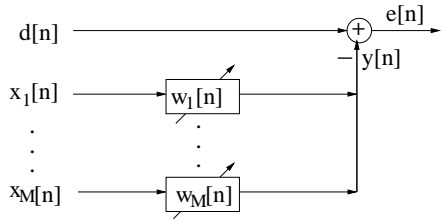}
   \caption{Structure of a reference signal based adaptive beamformer.
    \label{fig:multi_time_structure}}
\end{center}
\end{figure}
When a reference signal  $d[n]$ is available, adaptive beamforming can be implemented by the standard adaptive filter structure, as shown in Fig.~\ref{fig:multi_time_structure}, where $x_m[n]$, $m=1, 2, \cdots, M$ are the received quaternion-valued vector sensor signals,  $w_m[n]$, $m=1, 2, \cdots, M$ are the corresponding quaternion-valued coefficients, $y[n]$ is the beamformer output and $e[n]$ is the error signal.


\section{Simulation Results}\label{sec:simulations}

\subsection{Scenario one}
\begin{figure}[t]
\begin{center}
   \includegraphics[width=0.9\linewidth]{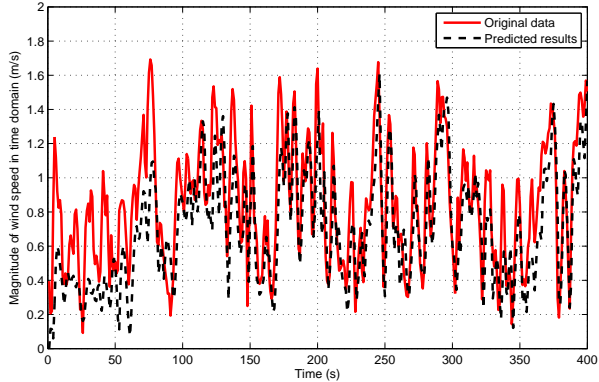}
   \caption{Prediction results using the QLMS algorithm.
    \label{fig:time_signal}}
\end{center}
\end{figure}
\begin{figure}
\begin{center}
   \includegraphics[width=0.9\linewidth]{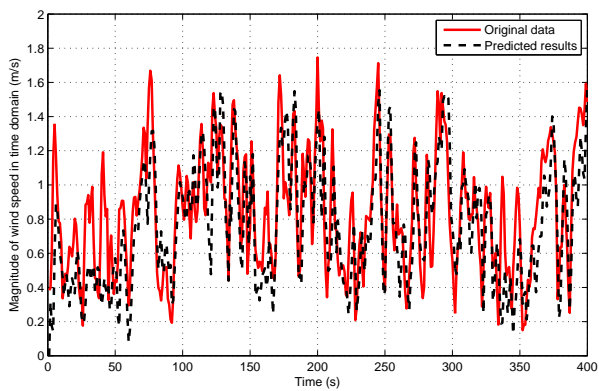}
   \caption{Prediction results using the AQLMS algorithm.
    \label{fig:aug_time_signal}}
\end{center}
\end{figure}
In this part, both the QLMS and the AQLMS algorithms are applied to the wind
data generated by CFD simulations with a sampling frequency of 1 Hz.
The parameters are as follows. The step size is $\mu=2.5\times10^{-4}$ and the adaptive filter length is $L=16$. The prediction step is 2. The adaptive weight vector is initialized as an all-zero vector. Fig. \ref{fig:time_signal} and Fig. \ref{fig:aug_time_signal} show the results for the QLMS and AQLMS algorithms, respectively. As we can see from the results, both algorithms can track the change of the wind speed signal effectively.

\subsection{Scenario two}
\begin{figure}
\begin{center}
   \includegraphics[width=0.9\linewidth]{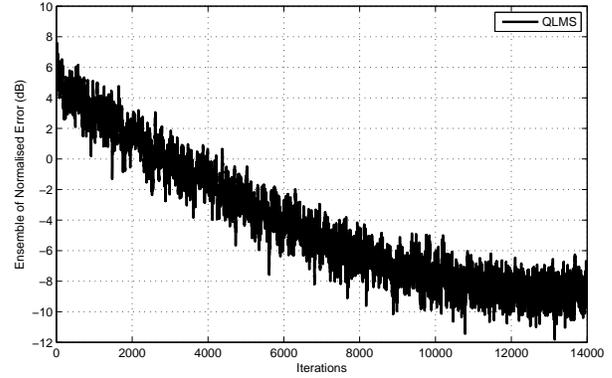}
   \caption{Learning curve using the  QLMS algorithm for adaptive beamforming.
    \label{fig:learning_curve_qlms}}
\end{center}
\end{figure}
Now we run simulations for the adaptive beamforming scenario. The vector sensor array with 10 crossed-dipoles and half-wavelength spacing is considered to obtain the output using the QLMS algorithm.  The stepsize  $\mu$ here is set to be $1\times10^{-6}$. A desired signal with 20 dB SNR impinges from the broadside and two interfering signals with the signal to interference ratio (SIR) of 0 dB arrive from $30^\circ$ and $-20^\circ$, respectively. All the signals have the same polarisation of $(\gamma, \eta)=(0,0)$. The learning curve averaged over 100 simulation runs is shown in  Fig.~\ref{fig:learning_curve_qlms}, where we can see the normalised error has reached about -10 dB, indicating an effective beamforming operation.

\section{Conclusion}\label{sec:conclusions}

In this paper, a general quaternion-valued gradient operator has been derived in detail, based on which  two adaptive algorithms were developed including the QLMS and the  AQLMS algorithms. These algorithms were applied to two different areas.   One  is to combine with the classic computational fluid dynamics (CFD) approach in wind profile prediction and the other one is to apply the result to the adaptive beamforming problem for vector sensor arrays. Simulation results have shown that the derived algorithms can work in different scenarios effectively, highlighting the importance and usefulness of the derived gradient operator. One important note is that  although there have been some derivations of this operator in literature with different level of details, this is the first time to give the most general form with a solid theoretical basis.

\section{Acknowledgements}
This work is partially funded by National Grid UK.



\end{document}